\documentclass[11pt,reqno]{article}

\usepackage[usenames]{color}
\usepackage{graphicx}
\usepackage{amscd}
\usepackage[english]{babel}
\usepackage{color}
\usepackage{fullpage}
\usepackage{psfig}
\usepackage{graphics,amsmath,amssymb}
\usepackage{amsthm}
\usepackage{amsfonts}
\usepackage{latexsym}
\usepackage{epsf}
\usepackage{float}
\usepackage{MnSymbol}

\usepackage[colorlinks=true,
linkcolor=webgreen,
filecolor=webbrown,
citecolor=webgreen]{hyperref}

\definecolor{webgreen}{rgb}{0,.5,0}
\definecolor{webbrown}{rgb}{.6,0,0}
		
\theoremstyle{plain}
\numberwithin{equation}{section}

 \newcommand{\seqnum}[1]{\href{http://oeis.org/#1}{\underline{#1}}} 
 
\begin{document}

\setcounter{page}{1}


\begin{center}

{\large\bf The {\em OEIS}: A Fingerprint File for Mathematics  } \\
\vspace*{+.2in}

N. J. A. Sloane\footnote{N. J. A. Sloane is President of The OEIS Foundation, Inc., and is a Visiting Scholar at
the Math. Dept., Rutgers University, New Brunswick, NJ. His email address is njasloane@gmail.com.} \\
\vspace*{+.1in}

\end{center}


\section{The OEIS}\label{SecI}
The {\em On-Line Encyclopedia of Integer Sequences} or {\em OEIS} \cite{OEIS} is a free website (\url{https://oeis.org}) containing information about 350,000 number sequences. You will probably first encounter it when trying to identify a sequence that has come up in your work. If your sequence is recognized, the response will  tell you the first 100 or sometimes  10,000 terms, give a definition, properties, formulas, references, links, computer programs, etc., as appropriate. If it is not recognized, the system will invite you to submit it if you think it is of general interest, so that the next person who looks it up will find it - and be grateful to you.

Ron Graham called the OEIS a ``fingerprint file for mathematics''. It has also been described as, ``pound for pound, one of the most useful mathematics sites on the web''. If you have been struggling with a sequence, and the OEIS tells you what it is, you will understand why it is so popular.  Of course, if it tells you that your problem was solved forty years ago, you may be unhappy, but it is better to find out now rather than later. The database has been around, in one form or another, for 57 years, so if your sequence is not yet included, there is a moderate chance it is new (of course this is not a proof).

\section{The OEIS as a source of problems}\label{Sec2Prob}
The entries are constantly being updated.
Every day we get a hundred or so submissions of new sequences, and another hundred 
comments on existing entries (new formulas, references, additional terms, etc.).

The new sequences  are often sent in by non-mathematicians,
and are a great source of problems to work on. You can see the current submissions at \url{https://oeis.org/draft}.
Often enough you will see a sequence that is so interesting you want to drop everything else and work on it.
Well, go ahead!  Many research papers have been born in this way. The ``Yellowstone Permutation''
\seqnum{A098550}
is one of my favorite examples: it was quite a challenge to prove that it {\em is}
a permutation of the positive integers, and the entry has a link to a paper \cite{Yellow}  that a group of us wrote
analysing it.

Some of the things you might work on after seeing an interesting sequence on the drafts stack are:
\begin{itemize}
\item Is the sequence well-defined?
\item Does it contain infinitely many terms?
\item How fast does it grow? The ``graph'' button in every entry is helpful here.
\item Is there a formula or generating function?
\item Or you may see a conjecture or question  that you think you can answer.
\item How would you program it to generate more terms?
\end{itemize}
Studying the ``drafts''  stack is an endless source of fun (and sleepless nights).

\section{The OEIS as a reference work}\label{SecRef}
Even when you know what the sequence is, you may still look it up, to find out what is presently known about it. 
We try to make sure the sequences are well-supplied with references, especially any recent articles.
The coverage is broad: besides the obvious fields like combinatorics, graph theory, 
group theory, number theory, computer science, recreational mathematics, there are
large numbers of sequences from physics and chemistry.

This makes the OEIS extremely useful as a reference work.
Earlier this year a senior group theorist told us he had been working on a problem for 20 years, but when he looked it up in the OEIS he found a reference that he was not aware of. 

Another great resource are the computer programs. An important sequence will have programs to generate it in C, Java, MAGMA, Maple, Mathematica, PARI, Python, Sage, etc. This is very helpful for numerical investigations.

\section{Submitting your sequence or comment}\label{Sec2Sub}
If you have an interesting number sequence that is not in the OEIS, you should definitely submit it.
Having a sequence in the OEIS is something you can be proud of.
You will be joining an enterprise that has been running
for almost 60 years, and to which over 12,000 people have already contributed.
Contributions come from almost every country, and we have been called one of the most successful 
international collaborations.

You must register before you can contribute, using your real name, 
the name you would use in a scientific publication 
 (see \href{https://oeis.org/wiki/Special:RequestAccount}{Request Account}).
  The OEIS is not a ``social media''. One of
 the reasons for  the success of the database is that we have high standards, all contributions are
 refereed by the editors, and accuracy is of primary importance.
 See \href{https://oeis.org/wiki/Overview_of_the_contribution_process}{here} for more about the submission process.
 The OEIS Wiki also has a page showing examples of
\href{https://oeis.org/wiki/Examples_of_what_not_to_submit}{What not to submit}!
 
\section{The OEIS Wiki}\label{SecW}
The OEIS Wiki (\url{https://oeis.org/wiki}) has a great deal of useful information for  users, especially in the ``Information'' section. There is a general index to the sequences, a style sheet, a Q\&A page, an FAQ page, pages called 
\href{https://oeis.org/wiki/How_to_add_a_comment,_more_terms,_or_a_b-file_(short_version)}{How to add a comment, more terms, or a b-file, ...},
\href{https://oeis.org/wiki/Instructions_For_General_Users}{Instructions for general users}, 
\href{https://oeis.org/wiki/The_multi-faceted_reach_of_the_OEIS}{The multi-faceted reach of the OEIS}, and so on.

\section{Citations of the OEIS}\label{SecCite}
An especially important part of the OEIS Wiki is the section
 (see \url{https://oeis.org/wiki/Works_Citing_OEIS}) that lists citations of the OEIS in
the literature or on the web. There are now about 10,000 citations, which often say things like ``This theorem would not have been discovered
without the help of the OEIS''.

You can help keep these pages up-to-date. If you come across a paper that mentions a number sequence,
check if the sequence is in the OEIS, add it if it isn't, and make sure the entry has a reference to the paper.
If the context seems new, consider adding a comment saying ``Arises in the spectral analysis of cobweb singularities, 
see  A. Spider et al. (2021)."  If you come across an article that references the OEIS, make sure the paper 
is listed in the ``Works  Citing the OEIS'' pages on the Wiki.  And if you happen to see a comment that ``This sequence is not in the OEIS", then add it at once.

When you write a paper yourself, using information from the OEIS, don't forget to mention us in your references list,
typically by saying 

The OEIS Foundation Inc., Entry A123456. Published electronically at https://oeis.org, 2021,

\noindent
and also mention it in any relevant OEIS entries.
Many authors ``forget'' to do this, but it will help your career by drawing attention to your paper.

\section{``Link Rot''}\label{SecL}

Don't get me started! The OEIS serves as a guide to the literature on a huge number of subjects,
and we have hundreds of thousands of links. And every day many of these links break.
System administrators feel they are not doing their job if they don't change all their URLs every
couple of years. Or, as happened last year, a major university will decide to delete all the faculty  home pages, with no warning. That broke several hundred of our links.
Pages on individual's web sites are the most fragile of all.
 
If you are a frequent user of the OEIS you will often run into this problem. You can help by locating a replacement
URL if the site has simply moved, or by adding a new link to a copy of the missing page on the wonderful Wayback Machine, run by the Internet Archive \cite{WBM}.  A better solution is to ask the author of the page 
for permission to put a local copy of the page on the OEIS server. We have a strong reputation, and hope to be around for a long time. Almost everyone we have asked has agreed.
But it is a time-consuming business, and you could help.

\section{Other topics}\label{SecOT}

\vspace*{+.2in}
\noindent{\bf ``Superseeker''.} This is a program that runs on our server, and tries very hard to identify a sequence.
It runs several programs that try to guess a formula or recurrence (including the powerful 
Salvy-Zimmermann program  {\em gfun} \cite{gfun}). It also transforms 
the sequence in a hundred ways and looks up the results in the OEIS,  hoping for a match.
To use it, send an email 
to  \href{mailto:superseeker@oeis.org}{\tt superseeker@oeis.org}.
with a blank subject, containing a single line of the form
   
   lookup 0 1 3 7 11 15 23 35 43 47
   
\noindent
(with no commas). 
 Since this uses many resources on our server, please use it sparingly.
   
   One of my current goals is to strengthen Superseeker. If you are interested in helping, please contact me.

   \vspace*{+.2in}
   \noindent{\bf The Sequence Fans Mailing List.}  Any registered user of the OEIS can join, and 
   messages go out to a few hundred fellow sequence-lovers.  This can be even more powerful than Superseeker if you are really desperate to identify a sequence. There is a link on the OEIS Wiki.

    \vspace*{+.4in}
   \noindent{\bf Triangles and arrays of numbers.} 
   The OEIS also includes triangles and arrays of numbers. Triangles are read by rows. Pascal's triangle
   becomes \seqnum{A007318}: $1, ~1, 1, ~1, 2, 1,$ $1, 3, 3, 1, ~1, 4, 6, 4, 1,$ $1, 5, 10, 10, 5, 1, \ldots$.
   Arrays are read by antidiagonals. The table of Nim-sums $m \oplus n$ \cite{BON},
$$
\begin{array}{ccccccc}
0 & 1 & 2 & 3 & 4 & 5 & \ldots \\
1 & 0 & 3 & 2 & 5 & 4 & \ldots  \\
2 & 3 & 0 & 1 & 6 & 7 & \ldots  \\
3 & 2 & 1 & 0 & 7 & 6 & \ldots  \\
4 & 5 & 6 & 7 & 0 & 1 & \ldots  \\
. & . & . & . & . & . & \ldots  \\
\end{array}
$$ 
becomes \seqnum{A003987}, $0, ~1, 1,~ 2, 0, 2,$ $ 3, 3, 3, 3, $ $ 4, 2, 0, 2, 4, \ldots$.

\vspace{0.2in}

The OEIS aims for a broad coverage of integer sequences arising in science, especially in mathematics (of course including several versions of the famous subway stops sequence, since they have been published on tests;  but not the number of pages in the $n$th issue of these {\em Notices} since 1954, because no one has mentioned it before). Our motto is accuracy and completeness, combined with good judgment.
We have a great many contributors, a great many entries, and we hope you will join us at \url{https:oeis.org}.


\bigskip
\hrule
\bigskip

\bigskip

\noindent 2020 Mathematics Subject Classification: 00Axx, 05-00, 11-00, 11Bxx, 97K10



\begin{thebibliography}{99}

\bibitem{Yellow} 
David L. Applegate, Hans Havermann, Bob Selcoe, Vladimir Shevelev, N. J. A. Sloane, and Reinhard Zumkeller, 
The Yellowstone Permutation, 
\emph{J. Integer Sequences}, \textbf{18:6} (2015), \#15.6.7.

\bibitem{BON} 
J. H. Conway and R. K. Guy,
\emph{The Book of Numbers},
Springer, 1996; p.~295.

\bibitem{WBM} 
Internet Archive,
\emph{The Wayback Machine},
\tt{https://archive.org/web},  \rm 2021.

\bibitem{OEIS} 
The OEIS Foundation Inc.,
\emph{The On-Line Encyclopedia of Integer Sequences},
\tt{https://oeis.org} \rm (2021).

\bibitem{gfun} 
Bruno Salvy and Paul Zimmermann, 
GFUN: a Maple package for the manipulation of generating and holonomic functions
in one variable, 
\emph{ACM Trans. Math. Software}, \textbf{20} (1994), 163--177.


\end{thebibliography}
\end{document}